\documentclass[12pt]{amsart}
\usepackage{amsmath,amsfonts,amssymb,amsthm}
\usepackage{pinlabel}
\usepackage{graphicx}
\usepackage{mathtools}
\usepackage[all]{xy}
\usepackage{hyperref}
\usepackage[usenames,dvipsnames]{color}
\usepackage{color}
\xyoption{dvips}

\newtheorem{theorem}{Theorem}

\newtheorem{proposition}[theorem]{Proposition}

\theoremstyle{remark}
\newtheorem{remark}[theorem]{Remark}

\theoremstyle{example}
\theoremstyle{definition}

\newcommand{\R}{\mathbb{R}}

\numberwithin{equation}{section}

\begin{document}

\title[Non-decomposable Lagrangian cobordisms between Legendrian knots]{Non-decomposable Lagrangian cobordisms between Legendrian knots}

\author{Roman Golovko}

\author{Daniel Kom\'{a}rek}

\begin{abstract}
For a given $g>0$, we construct a family of non-decomposable Lagrangian cobordisms of genus $g$ between (stabilized) Legendrian knots in the standard contact three-sphere. The main technique we use to obstruct decomposability is based on Livingston's estimates from \cite{Critpointsknotcobordisms}.
\end{abstract}

\address{Faculty of Mathematics and Physics, Charles University, Sokolovsk\'{a} 49/83, Prague 8, 186 00,  Czech Republic}
\email{golovko@karlin.mff.cuni.cz}

\address{Faculty of Mathematics and Physics, Charles University, Sokolovsk\'{a} 49/83, Prague 8, 186 00 , Czech Republic}
\email{daniel.komarek@matfyz.cuni.cz}

\date{\today}
\thanks{}


\maketitle

\section{Introduction and main results}

The standard way to construct Lagrangian cobordisms between Legendrian knots in the standard contact $S^3_{std}$ is by concatenating elementary 
cobordisms which are traces of three moves in the front projection: Legendrian isotopy,
isolated standard unknot birth (induces Lagrangian 0-handle attachment), and Legendrian surgery, which is called a pinch move in reverse (it induces Lagrangian 1-handle attachment). See Figure 1 for the schematic description of the latter two pieces. Cobordisms obtained this way are called decomposable and have been first constructed by Chantraine in \cite{SomeNonColSl} and Ekholm-Honda-Kalman in \cite{LegKnotsLagrCob}.  

Decomposable Lagrangian cobordisms can be seen as symplectic analogues of ribbon cobordisms in smooth topology. 
Recall that for a given two knots $K_0,K_1\subset S^3$ a ribbon cobordism $\Sigma$ from $K_0$ to $K_1$ is a smooth cobordism $\Sigma\subset S^3\times[0,1]$ from $K_0$ to $K_1$ such that the projection $S^3\times [0,1]\to [0,1]$ restricts to a Morse function on $\Sigma$ which does not have index two critical points.
Note that all decomposable Lagrangian cobordisms between Legendrian
knots in $S^3_{std}$ are ribbon.

\begin{figure}[t]
\begin{center}
\vspace{3mm}
\labellist
\pinlabel $t$ at 220 585
\endlabellist
\includegraphics[width=220px]{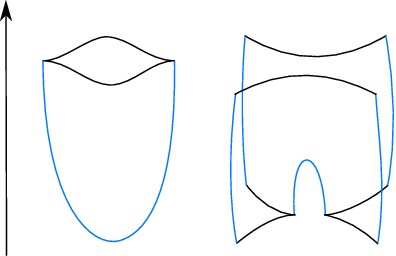}
\caption{Left: Lagrangian filling of the $tb=-1$ unknot; Right: Lagrangian pair-of-pants cobordism obtained from the Lagrangian $1$-handle attachment.}
\label{fig:decomposable}
\end{center}
\end{figure}

These days, non-decomposable Lagrangian cobordisms are proven to exist in the following two forms: 
\begin{itemize}
\item non-decomposable Lagrangian caps with stabilized
Legendrian ends that
have been discussed in the work of Lin \cite{LinExLagrCaps}; 
\item non-decomposable Lagrangian concordances between stabilized Legendrian knots \cite{NonregStabilizedConc} and between Lagrangian fillable Legendrian knots \cite{Nonregfillable} that have been constructed by Dimitroglou Rizell and the first author.
\end{itemize}
In this note we construct non-decomposable Lagrangian cobordisms of an arbitrary non-trivial genus between Legendrian knots in $S^3$. The method to determine non-decomposability that we use in this paper is different from the ones we used with Dimitroglou Rizell in \cite{NonregStabilizedConc}, \cite{Nonregfillable}. Recall that in \cite{NonregStabilizedConc} we use the obstruction of non-decomposability based on the variation of the result of  Cornwell-Ng-Sivek \cite[Theorem 3.2]{ObstrLagrConc}, and in \cite{Nonregfillable} we relied on the result of Agol which says that ribbon (stably homotopy ribbon) Lagrangian concordances of knots form a partial order, see \cite[Theorem 1.1]{RibbonConcPartialOrdering}. In this paper, we rely on the obstruction of ribbon cobordisms due to Livingston \cite[Corollary 5.4]{Critpointsknotcobordisms} and its partial adaptation due to Borodzik and Tru\"ol \cite[Theorem 2.5]{BorodzikTruol25}.

The main result of this paper is the following:
\begin{theorem}
\label{mainthgenusg}
Let  $\Lambda$ be a Legendrian knot in $(S^3, \xi_{st})$ such that
\begin{itemize}
\item[(1)] there is an odd prime $p$ that divides the determinant of $\Lambda$, 
\item[(2)] there is a decomposable  Lagrangian concordance $C$ from the $tb=-1$ unknot $U$ to  $\Lambda$. 
\end{itemize}
Then for a given $g > 0$, and for all integers $n > 2g$, there exists a non-decomposable
Lagrangian cobordism $L_g$ of genus $g$, from a Legendrian knot $\Lambda_{-}$ obtained from $\Lambda^n:=\Lambda\#\dots\#\Lambda$
by sufficiently many stabilizations, to a Legendrian knot $\Lambda_+$ obtained from $U^n:=U\#\dots\#U$
by sufficiently many stabilizations.
\end{theorem}

\begin{remark}
Note that $\det(U)=1$, and hence condition $(1)$ in the formulation of Theorem \ref{mainthgenusg} implies that $\Lambda$ is not smoothly isotopic to $U$.
\end{remark}

\section{Applicability of Theorem \ref{mainthgenusg}}
Here we describe a few series of examples to which we can apply Theorem \ref{mainthgenusg}.
\subsubsection*{Examples A} 
\begin{figure}[h]
\begin{center}
\vspace{3mm}
\labellist
\pinlabel $t$ at 120 705
\pinlabel $k-3$ at 420 670
\endlabellist
\includegraphics[width=200px]{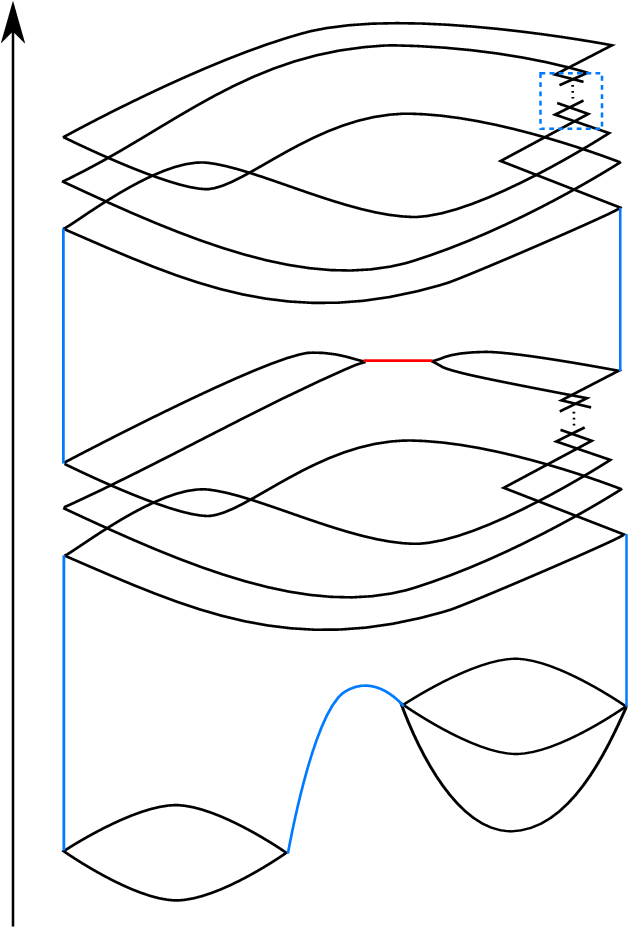}
\caption{Decomposable Lagrangian concordance $C_{k}$ from the $tb=-1$ Legendrian unknot $U$  to the Legendrian representative $\Lambda_{k}$ of the pretzel knot $P(3,-3,-k)$, where $k\geq 3$, induced by the ambient surgery along the red arc.   The front projection of $\Lambda_{k}$ has $k-3$ crossings in the blue box.} 
\label{fig:concordance_pretzel}
\end{center}
\end{figure}

Consider the sequence of Legendrian pretzel knots $P(3,-3,-k)$, where $k\geq 3$, from \cite{DoublysliceandObstr,NonregStabilizedConc}, we call them $\Lambda_{k}$. 
There is a decomposable Lagrangian concordance $C_{k}$ from $U$ to  $\Lambda_{k}$,  see Figure \ref{fig:concordance_pretzel}, this concordance has been discussed by Chantraine and Legout in \cite{DoublysliceandObstr} 
and by Dimitroglou Rizell and the first author in \cite{NonregStabilizedConc}. 
Following the computation of the determinant for pretzel knots from \cite{BaeLeeAlexanderPretzel} we observe that 
$\det(P(3,-3,-k))=9$ and the odd prime $p$ divisor of  $\det(P(3,-3,-k))$ that is required  by Theorem \ref{mainthgenusg} can be taken to be $p=3$. Hence, Theorem \ref{mainthgenusg} 
can be applied to $C_k$.

\subsubsection*{Examples B}
In order to describe the following family we first recall a few properties of determinant and Alexander polynomial.
\begin{remark}
\label{behaviorAlexanderdet}
Determinant and Alexander polynomial satisfy the following properties:
\begin{itemize}
\item[1)] For a given knot $K$, $|\Delta_{K}(-1)|=\det(K)$.
\item[2)] For given two knots $K_1$ and $K_2$, the Alexander polynomial of the connected sum  $K_1\# K_2$ is given by 
$\Delta_{K_1\# K_2}(t) = \Delta_{K_1}(t)\Delta_{K_2}(t)$, which implies that 
$$\det(K_1\# K_2) = |\Delta_{K_1\# K_2}(-1)| = |\Delta_{K_1}(-1)||\Delta_{K_2}(-1)| = \det(K_1)det(K_2).$$
\item[3)] For a knot $K$, $\det(K)\neq 0$. This, in particular, can be justified by the fact that $\det(K)$ can be seen as the order of the finite abelian group $H_1(\Sigma_2(K);\mathbb Z)$, 
where $\Sigma_2(K)$ denotes the 2-fold branched cover of $S^3$ over $K$. 
\end{itemize}
\end{remark}
Now we take a Legendrian knot $\Lambda$ such that there is a decomposable Lagrangian concordance $C'$ from $U$ to $\Lambda$. Then
we take the connected sum of decomposable Lagrangian concordances $C_{k}\# C'$ from $U\# U\simeq U$ to $\Lambda_{k}\# \Lambda$. Such an operation can be performed following the discussion in \cite[Section 5]{LegSatDecomLagrCob}. Observe that $C_{k}\# C'$ is decomposable by \cite[Corollary 5.1, Remark 5.3]{LegSatDecomLagrCob}. 
Then note that from the third part of  Remark  \ref{behaviorAlexanderdet} it follows that $\det(\Lambda)\neq 0$. 
From this, together with the second part of Remark \ref{behaviorAlexanderdet}, it follows that 
$p=3 |\det(\Lambda_{k}\# \Lambda)$. So, we can apply Theorem \ref{mainthgenusg} to the concordance $C_{k}\# C'$.
\begin{remark}
For example, we can take $C'$ to be a decomposable Lagrangian concordance from $U$ to the Legendrian representative of $(2,n)$-torus knot from \cite{ExactLagrfillingsof2n}.
\end{remark}
\subsubsection*{Examples C} Here we rely on the decomposable realization of ribbon cobordisms developed by Etnyre and Leverson in \cite{EtnyreLeversonRealization}. We can simply take a smooth knot $K$ satisfying  that $p | det(K)$ for some $p>1$ and such that there is a ribbon concordance from an unknot to $K$. Then following \cite[Theorem 1.2]{EtnyreLeversonRealization} we realize this concordance by a decomposable Lagrangian concordance $C$ from the Legendrian representative of an unknot (in general, stabilized) to    the  Legendrian representative of $K$ (in general, stabilized). Obviously, Theorem \ref{mainthgenusg} can be applied to $C$.

\section{Proof of Theorem \ref{mainthgenusg}}
We start the proof by considering  the connected sum of $n$ copies of concordance $C$ and denote it by $C^n:=C\#\ \dots \# C$. 
Since $C$ is a decomposable Lagrangian concordance, from \cite[Corollary 5.1, Remark 5.3]{LegSatDecomLagrCob} it follows that $C^n$ is a decomposable Lagrangian concordance from the Legendrian knot $U^n=U\#\dots \#U$, which is smoothly an unknot, to the Legendrian knot $\Lambda^n=\Lambda\# \dots \# \Lambda$. We know that there is a prime number $p$ that divides  $\det(\Lambda)$. From the second part of Remark \ref{behaviorAlexanderdet} we observe that $\det(\Lambda^n)=\det(\Lambda)^n$, and hence 
$p$ that divides  $\det(\Lambda)$ also divides $\det(\Lambda^n)$ and does not divide $\det(U^n)=\det(U)=1$.

Then we argue the same way as in \cite[Section 2]{Nonregfillable} and \cite[Section 5]{NonregStabilizedConc}. First, we 
revert the concordance $(C^n)^{-1}=\{(-t,x)\in \R_t\times \R^3_x\ | \ (t,x)\in C^n\}.$
 Since $C^n$ is decomposable, observe that the tangent planes of $(C^n)^{-1}$ are homotopic to totally-real tangent planes. Using this, we apply the $h$-principle for totally real embeddings and see that $(C^n)^{-1}$  admits a compactly supported smooth isotopy to a concordance that is totally real for any choice of cylindrical almost complex structure, we denote it by $T^n$.
Finally, we modify $T^n$ using the approximation result of Dimitroglou Rizell \cite{LagrApproxTotReal}
and get the Lagrangian concordance $S(T^n)$ from $S(\Lambda^n)$ to $S(U^n)$. Here $S(U^n)$ and $S(\Lambda^n)$ are multiply stabilized versions of $U^n$ and $\Lambda^n$, respectively.

\begin{figure}[h]
\begin{center}
\vspace{3mm}
\labellist
\pinlabel $t$ at 235 590
\pinlabel $\Lambda'$ at 430 550
\pinlabel $(S_+S_-)(\Lambda')$ at 450 400
\endlabellist
\includegraphics[width=190px]{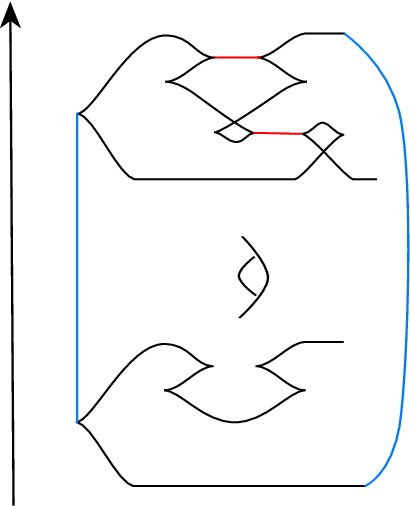}
\caption{Lagrangian cobordism of genus one from $(S_+S_-)(\Lambda')$ to $\Lambda'$ obtained by making an ambient surgery along two red arcs.}
\label{fig:genusonedoublestab}
\end{center}
\end{figure}

\subsection{Lagrangian cobordisms constructed by Sabloff, Vela-Vick and Wong}
\label{highgenusfromstab}
Sabloff, Vela-Vick and Wong in \cite{SabloffVelaVickWongUpperBound} constructed the following cobordism that will be useful for us. 
Given a Legendrian knot $\Lambda'\subset S^3_{std}$, let $(S_+S_-)(\Lambda')$ be the Legendrian knot obtained by first applying the negative stabilization to $\Lambda'$ followed by the application of the positive stabilization. Sabloff, Vela-Vick and Wong in \cite[Lemma 3.2]{SabloffVelaVickWongUpperBound} showed that there is a  Lagrangian cobordism 
of genus one from  $(S_+S_-)(\Lambda')$ to $\Lambda'$, the schematic picture of it appears in Figure \ref{fig:genusonedoublestab}. 

We apply this construction inductively to $\Lambda'=S(\Lambda^n)$.
In other words, given $g > 0$,
we concatenate $g$ copies of the genus-$1$ Lagrangian cobordism from \cite[Lemma 3.2]{SabloffVelaVickWongUpperBound} 
to obtain a genus-$g$ Lagrangian cobordism from $(S_+S_-)^g(S(\Lambda^n))$ to $S(\Lambda^n)$. 
Note that the genus is additive because the cobordisms that we are concatenating have connected boundaries.

We then concatenate the obtained cobordism and $S(T^n)$ and get a Lagrangian cobordism $L_g$ of genus $g$ from $\Lambda_-=(S_+S_-)^g(S(\Lambda^n))$ to $\Lambda_+=S(U^n)$.
We claim that $L_g$ is a non-ribbon cobordism. In order to show it, we need to rely on the estimates of Livingston from \cite{Critpointsknotcobordisms}.

\subsection{Livingston's estimates and the results of Borodzik--Tru\"ol}
For a given two knots $K_0$ and $K_1$ in $S^3$ and a connected Morse cobordism  $\Sigma\subset S^3\times [0,1]$ from $K_0$ to $K_1$ (i.e. $\Sigma$ is a connected cobordism satisfying that the projection of $\Sigma$ to $[0,1]$ is a Morse function), Livingston in \cite{Critpointsknotcobordisms} investigates the behaviour of quadruples $(c_0,c_1,c_2,g)$, where $c_i=c_i(\Sigma)$ denotes the number of critical points of index $i$ on $\Sigma$ and $g$ is a genus of $\Sigma$. More precisely, Livingston in  \cite{Critpointsknotcobordisms} proved a few obstructions to the existence of $\Sigma$ in the form of estimates involving the homology of the $n$-fold branched covers $\Sigma_n(K_i)$ of $S^3$ over $K_i$, $i=0,1$.

 In this note we rely on the partial case of Livingston's estimates that uses the homology of the two-fold branched covers $\Sigma_2(K_0)$ and $\Sigma_2(K_1)$.
\begin{theorem}[\cite{Critpointsknotcobordisms}]
If $\Sigma$ is a connected Morse cobordism from $K_0$ to $K_1$, then for an arbitrary odd prime $p$ the following inequality holds
\begin{align}
\label{livingstonsestimatec2}
c_2(\Sigma)\geq \frac{\beta_1(\Sigma_2(K_0); \mathbb F_p)-\beta_1(\Sigma_2(K_1); \mathbb F_p)}{2}-g(\Sigma).
\end{align}
Here $\beta_1(\Sigma_2(K_i); \mathbb F_p)$ is the dimension of $H_1(\Sigma_2(K_i); \mathbb F_p)$ as
an $\mathbb F_p$-vector space, $i=0,1$.
\end{theorem}

Note that our handle convention is opposite to that of Livingston in \cite{Critpointsknotcobordisms}: Livingston considers cobordisms from $K_1$ to $K_0$, whereas we consider cobordisms from $K_0$ to $K_1$, and the handle indices are reversed.

Besides that, we need the following statement. It is a standard exercise, and we refer to the work of Borodzik--Tru\"ol \cite{BorodzikTruol25} for a proof.

\begin{proposition}[\cite{BorodzikTruol25}]
\label{basicestimates}
Given a knot $K\subset S^3$ and a prime $p$.   
\itemize
\item[(1)]
If $p$ divides $\det(K)$, then $\dim_{\mathbb F_p}H_1(\Sigma_2(K^n=K\#\dots \# K); \mathbb F_p)\geq n$.
\item[(2)]
If $p$ does not divide $\det(K)$, then $H_1(\Sigma_2(K); \mathbb F_p)=0$.
\end{proposition}

\subsection{Application of Livingston's estimates to obstruct decomposability}
Now we take the Lagrangian cobordism $L_g$ of genus $g$ from $(S_+S_-)^g(S(\Lambda^n))$ to $S(U^n)$  that we constructed in Section \ref{highgenusfromstab}.
From the second part of Remark \ref{behaviorAlexanderdet} we observe that 
$$\det((S_+S_-)^g(S(\Lambda^n)))=\det(\Lambda)^n\quad \mbox{and}\quad  \det(S(U^n))=1.$$ Besides that, $(S_+S_-)^g(S(\Lambda^n))$ is smoothly isotopic to $\Lambda^n$.
Hence, using the fact that $p| \det(\Lambda)$, from the first part of  Proposition \ref{basicestimates} it follows that 
\begin{align}
\label{nontton}
\dim_{\mathbb F_p}H_1(\Sigma_2((S_+S_-)^g(S(\Lambda^n))); \mathbb F_p)=\dim_{\mathbb F_p}H_1(\Sigma_2(\Lambda^n); \mathbb F_p)\geq n.
\end{align}
Besides that, the second part of Proposition \ref{basicestimates} implies that 
\begin{align}
\label{trton}
H_1(\Sigma_2(S(U^n)); \mathbb F_p)=H_1(\Sigma_2(U); \mathbb F_p)=0,\quad \mbox{and hence}\quad \dim_{\mathbb F_p} H_1(\Sigma_2(U); \mathbb F_p)=0. 
\end{align}
Finally, combining Formulas \ref{livingstonsestimatec2}, \ref{nontton} and \ref{trton}, we get that 
$c_2(L_g)\geq \frac{n}{2}-g$. Therefore, for any fixed $g>0$, since $n>2g$, we see that $c_2(L_g)>0$.  So, we see that $L_g$ is a non-ribbon Lagrangian cobordism of genus $g>0$, and hence $L_g$ is a non-decomposable Lagrangian cobordism of genus $g>0$. This finishes the proof.

\begin{remark}
Note that the ends of Lagrangian cobordisms that we get in Theorem \ref{mainthgenusg} are stabilized. One could try to extend Theorem \ref{mainthgenusg} to Lagrangian cobordisms with fillable Legendrian ends as Dimitroglou Rizell and the first author did for Lagrangian concordances in \cite{Nonregfillable}. In \cite{Nonregfillable}, Dimitroglou Rizell and the first author used Legendrian Whitehead doubling to achieve this goal. Unfortunately, it is difficult to mimic the application of Legendrian Whitehead doubling in these settings. The reason for that is that in the proof of Theorem \ref{mainthgenusg} we strongly rely on the first part of Proposition \ref{basicestimates} which allows to control the size of $\dim_{\mathbb F_p}H_1(\Sigma_2(\cdot); \mathbb F_p)$. Legendrian Whitehead doubling operation unfortunately does not respect (commute with) the connected sum operation. This makes the control of 
$\dim_{\mathbb F_p}H_1(\Sigma_2(\cdot); \mathbb F_p)$ quite difficult.
\end{remark}

\section*{Acknowledgements}
The authors would like to thank Georgios Dimitroglou Rizell and Joshua Sabloff for the very helpful discussions. 
The first author was supported by GA\v{C}R Lead Agency Grant 26-20231L. The second author was supported by the grants GAUK 134125 and PRIMUS/24/SCI/009.


\begin{thebibliography}{}


\bibitem{RibbonConcPartialOrdering}
I.~Agol, {\em Ribbon concordance is a partial ordering}, Commun. Am. Math. Soc. 2 (2022), 374--379.

\bibitem{BaeLeeAlexanderPretzel}
Y.~Bae and I.~Lee, {\em On Alexander polynomials of pretzel links},  Kyungpook Mathematical Journal, 60(2):239--253, June 2020.

\bibitem{BorodzikTruol25}
M.~Borodzik and P.~Tru\"ol, {\em Non-complex cobordisms between quasipositive knots}, J. Math. Pures Appl, Volume 207
(2026), 103842.





\bibitem{SomeNonColSl}
B.~Chantraine, {\em Some non-collarable slices of Lagrangian surfaces}, Bull. Lond.
Math. Soc. 44 (2012), no. 5, 981--987.




\bibitem{DoublysliceandObstr}
B.~Chantraine and N.~Legout, {\em Doubly slice knots and obstruction to Lagrangian concordance},
Comptes Rendus. Mathématique, Volume 361 (2023), 1605--1609.



\bibitem{ObstrLagrConc}
 C.~Cornwell, L.~Ng, and S.~Sivek, {\em Obstructions to Lagrangian
concordance}, Algebr. Geom. Topol. 16 (2016), no. 2, 797--824.

\bibitem{LagrApproxTotReal}
G.~Dimitroglou Rizell, {\em Lagrangian approximation of totally real concordances}, Proc. London Math. Soc. (3) 2025;130:e70042. 

\bibitem{NonregStabilizedConc}
G.~Dimitroglou Rizell and R.~Golovko, {\em Instability of Legendrian knottedness, and non-regular Lagrangian concordances of knots}, preprint 2024. Available at arXiv:2409.00290.

\bibitem{Nonregfillable}
G.~Dimitroglou Rizell and R.~Golovko, {\em Non-regular Lagrangian concordances between Lagrangian fillable Legendrian knots}, preprint 2025. Available at arXiv:2509.13594.

\bibitem{LegKnotsLagrCob}
T.~Ekholm, K.~Honda, and T.~Kálmán, {\em Legendrian knots and exact Lagrangian
cobordisms}, J. Eur. Math. Soc. (JEMS) 18 (2012), 2627–2689.



\bibitem{EtnyreLeversonRealization}
J.~Etnyre and C.~Leverson, {\em Lagrangian realizations of ribbon cobordisms}, preprint 2024. Available at arXiv:2410.06305.




\bibitem{LegSatDecomLagrCob}
R.~Guadagni, J.~Sabloff, and M.~Yacavone, {\em Legendrian satellites and decomposable cobordisms}, Journal of Knot Theory and Its Ramifications, Vol. 31, No. 13, 2250071 (2022).



\bibitem{LinExLagrCaps}
F.~Lin, {\em Exact Lagrangian caps of Legendrian knots}, J. Symplectic Geom. 14(1):269–295 (2016).

\bibitem{Critpointsknotcobordisms}
C.~Livingston, {\em Critical point counts in knot cobordisms: abelian and metacyclic invariants}, Trans. Amer. Math. Soc. Ser. B, 10:171--200, 2023.




\bibitem{ExactLagrfillingsof2n}
Y.~Pan, {\em Exact Lagrangian fillings of Legendrian (2,n) torus links}, Pacific Journal of Mathematics, 289-2 (2017), 417--441.

\bibitem{SabloffVelaVickWongUpperBound}
J.~Sabloff, D.~Vela-Vick, and C.~Wong, {\em Upper bounds for the Lagrangian cobordism
relation on Legendrian links}, . Algebr. Geom. Topol. 24 (2024), no. 8, 4237--4263.

\end{thebibliography}
\end{document}